\documentclass[a4, 10pt]{amsart}
\usepackage{amssymb}
\usepackage{amstext}
\usepackage{amsmath}
\usepackage{amscd}
\usepackage{latexsym}
\usepackage{amsfonts}
\usepackage{graphicx}

\theoremstyle{plain}
\newtheorem{thm}{Theorem}[section]
\newtheorem*{thm*}{Theorem}
\newtheorem*{cor*}{Corollary}
\newtheorem*{prop*}{Proposition}

\newtheorem*{thma}{Theorem A}
\newtheorem*{thmb}{Theorem B}

\newtheorem{prop}[thm]{Proposition}
\newtheorem{lem}[thm]{Lemma}

\newtheorem{cor}[thm]{Corollary}

\newtheorem*{claim*}{Claim}

\theoremstyle{definition}
\newtheorem{defn}[thm]{Definition}
\newtheorem{ex}[thm]{Example}
\newtheorem{rem}[thm]{Remark}
\newtheorem*{rem*}{Remark}

\newtheorem*{conj*}{Conjecture}

\theoremstyle{remark}
\newtheorem*{pf}{{\sl Proof}}

\newtheorem*{tpf}{{\sl Proof of Theorem \ref{main}}}

\newtheorem*{ac}{{\sc Acknowledgments}}

\numberwithin{equation}{thm}
\def\Hom{\mathrm{Hom}}

\def\Ext{\mathrm{Ext}}

\def\Mod{\mathrm{Mod}}
\def\mod{\mathrm{mod}}

\def\Coker{\mathrm{Coker}}

\def\Im{\mathrm{Im}}

\def\rank{\mathrm{rank}}

\def\p{\mathfrak p}
\def\q{\mathfrak q}

\def\Z{\Bbb Z}

\def\D{{\mathcal D}}

\def\Supp{\mathrm{Supp}}
\def\Ann{\mathrm{Ann}}
\def\Ass{\mathrm{Ass}}

\def\Min{\mathrm{Min}}
\def\Max{\mathrm{Max}}

\def\height{\mathrm{ht}}
\def\grade{\mathrm{grade}}

\def\Spec{\mathrm{Spec}}

\def\A{{\mathcal A}}

\def\E{{\mathcal E}}

\def\X{{\mathcal X}}
\def\Y{{\mathcal Y}}
\def\W{{\mathcal W}}
\def\M{{\mathcal M}}
\def\L{{\mathcal L}}
\def\sss{{\mathcal S}}

\def\coh{\mathrm{coh}}
\def\Serre{\mathrm{Serre}}

\def\T{{\mathcal T}}
\def\zzz{{\mathcal Z}}
\def\kdim{\mathrm{dim}}

\def\lsp{{\mathcal L}_{\sf spcl}}
\def\ltors{{\mathcal L}_{\sf tors}}
\def\lserre{{\mathcal L}_{\sf Serre}}
\def\lcoh{{\mathcal L}_{\sf coh}}
\def\lthick{{\mathcal L}_{\sf thick}}
\def\lse{{\mathcal L}_{\sf subext}}

\def\dperf{\D_\mathrm{perf}}

\def\xx{\text{\boldmath $x$}}
\def\aa{\text{\boldmath $a$}}

\begin{document}

\title[Classifying subcategories of modules]{Classifying subcategories of modules over a commutative noetherian ring}
\author{Ryo Takahashi}
\address{Department of Mathematical Sciences, Faculty of Science, Shinshu University, 3-1-1 Asahi, Matsumoto, Nagano 390-8621, Japan}
\email{takahasi@math.shinshu-u.ac.jp}
\begin{abstract}
Let $R$ be a quotient ring of a commutative coherent regular ring by a finitely generated ideal.
Hovey gave a bijection between the set of coherent subcategories of the category of finitely presented $R$-modules and the set of thick subcategories of the derived category of perfect $R$-complexes.
Using this isomorphism, he proved that every coherent subcategory of finitely presented $R$-modules is a Serre subcategory.
In this paper, it is proved that this holds whenever $R$ is a commutative noetherian ring.
This paper also yields a module version of the bijection between the set of localizing subcategories of the derived category of $R$-modules and the set of subsets of $\Spec\,R$ which was given by Neeman.
\end{abstract}
\thanks{Part of this work was done while the author visited the University of Nebraska-Lincoln, partly supported by NSF grant DMS 0201904.}
\maketitle
\section{Introduction}

Around 1990, Hopkins \cite{Hopkins} and Neeman \cite{Neeman} gave a classification theorem of the thick subcategories of the derived category of perfect complexes (i.e. finite complexes of finitely generated projective modules) over a commutative noetherian ring in terms of the ring spectrum.
After that, Thomason \cite{Thomason} generalized this classification theorem to quasi-compact and quasi-separated schemes, in particular, to arbitrary commutative rings.
Let $\dperf(R)$ denote the derived category of perfect complexes over a commutative ring $R$.
The classification theorem (for commutative rings) can be stated as follows.

\begin{thm*}[Hopkins-Neeman-Thomason]
Let $R$ be a commutative ring.
Then there is an isomorphism
\begin{align*}
& \big\{\,\text{thick subcategories of }\dperf(R)\,\big\}\\
\,\cong\, & \big\{\,\text{complements of intersections of quasi-compact open subsets of }\Spec\,R\,\big\}
\end{align*}
of lattices.
\end{thm*}

Here we recall the definitions of several subcategories of an abelian category.
A coherent subcategory is defined to be a full subcategory which is closed under kernels, cokernels and extensions.
A Serre subcategory is defined to be a coherent subcategory which is closed under subobjects.
A torsion class is defined to be a Serre subcategory which is closed under arbitrary direct sums.
Let $R$ be a commutative ring.
We denote by $\Mod\,R$ the category of $R$-modules and by $\mod\,R$ the full subcategory of finitely presented $R$-modules.
If $R$ is noetherian, then the lattice of Serre subcategories of $\mod\,R$, the lattice of torsion classes of $\Mod\,R$, and the lattice of subsets of $\Spec\,R$ which are closed under specialization are isomorphic to each other.
Taking advantage of the Hopkins-Neeman-Thomason theorem, Garkusha and Prest \cite{GP,GP2} gave the following result very recently.
(A torsion class $\X$ of $\M=\Mod\,R$ is said to be of finite type if the inclusion functor $\M/\X\to\M$, where $\M/\X$ denotes the quotient category, preserves arbitrary direct sums.)

\begin{thm*}[Garkusha-Prest]
Let $R$ be a commutative ring.
Then the following hold.
\begin{enumerate}
\item[{\rm (1)}]
One has lattice isomorphisms
\begin{align*}
& \big\{\,\text{thick subcategories of }\dperf(R)\,\big\}\\
\,\cong\, & \big\{\,\text{torsion classes of finite type of }\Mod\,R\,\big\}\\
\,\cong\, & \big\{\,\text{complements of intersections of quasi-compact open subsets of }\Spec\,R\,\big\}.
\end{align*}
\item[{\rm (2)}]
Suppose that $R$ is coherent.
Then one has lattice isomorphisms
\begin{align*}
& \big\{\,\text{thick subcategories of }\dperf(R)\,\big\}\\
\,\cong\, & \big\{\,\text{Serre subcategories of }\mod\,R\,\big\}\\
\,\cong\, & \big\{\,\text{complements of intersections of quasi-compact open subsets of }\Spec\,R\,\big\}.
\end{align*}
\end{enumerate}
\end{thm*}

Also by using the Hopkins-Neeman-Thomason theorem, Hovey \cite{Hovey} proved the following classification theorem of coherent subcategories.

\begin{thm*}[Hovey]
Let $R$ be a quotient ring of a commutative coherent regular ring by a finitely generated ideal.
Then the following hold.
\begin{enumerate}
\item[{\rm (1)}]
One has a lattice isomorphism
\begin{align*}
& \big\{\,\text{thick subcategories of }\dperf(R)\,\big\}\\
\,\cong\, & \big\{\,\text{coherent subcategories of }\mod\,R\,\big\}.
\end{align*}
\item[{\rm (2)}]
Every coherent subcategory of $\mod\,R$ is a Serre subcategory.
\end{enumerate}
\end{thm*}

As Hovey pointed it out as an interesting fact, there was no direct proof of the second assertion of the above theorem; it could not be proved without resorting to the rather difficult classification of thick subcategories of the derived category, namely the Hopkins-Neeman-Thomason theorem.
Hovey conjectures that the isomorphism stated in the above theorem always holds for commutative coherent rings.
(Recall that a commutative ring is called coherent if every finitely generated ideal is finitely presented.)

\begin{conj*}[Hovey]
Let $R$ be a commutative coherent ring.
Then one has a lattice isomorphism
\begin{align*}
& \big\{\,\text{thick subcategories of }\dperf(R)\,\big\}\\
\,\cong\, & \big\{\,\text{coherent subcategories of }\mod\,R\,\big\}.
\end{align*}
\end{conj*}

One of the main purposes of this paper is to prove that this conjecture is true if $R$ is noetherian.
First of all, we will directly prove that every coherent subcategory of $\mod\,R$ is Serre; in the proof, we will not apply the Hopkins-Neeman-Thomason theorem.
Actually, we will not use the notion of a derived category in the proof.
Furthermore, the proof we will give is much simpler than Hovey's.
After that, we will prove that if $R$ is noetherian, then all the lattices appearing in the above part are isomorphic to each other.
Our first main theorem is the following.

\begin{thma}
Let $R$ be a commutative noetherian ring.
Then every coherent subcategory of $\mod\,R$ is a Serre subcategory, and one has the following isomorphisms of lattices:
\begin{align*}
& \big\{\,\text{thick subcategories of }\dperf(R)\,\big\}\\
\,\cong\, & \big\{\,\text{coherent subcategories of }\mod\,R\,\big\}\\
\,=\, & \big\{\,\text{Serre subcategories of }\mod\,R\,\big\}\\
\,\cong\, & \big\{\,\text{torsion classes of }\Mod\,R\,\big\}\\
\,\cong\, & \big\{\,\text{subsets of }\Spec\,R\text{ closed under specialization}\,\big\}\\
\,=\, & \big\{\,\text{complements of intersections of quasi-compact open subsets of }\Spec\,R\,\big\}.
\end{align*}
\end{thma}

On the other hand, Neeman \cite{Neeman} showed the following theorem.

\begin{thm*}[Neeman]
Let $R$ be a commutative noetherian ring.
Then one has an isomorphism
\begin{align*}
& \big\{\,\text{localizing subcategories of }\D(R)\,\big\}\\
\,\cong\, & \big\{\,\text{subsets of }\Spec\,R\,\big\}\\
\intertext{of lattices.
Moreover, this induces an isomorphism}
& \big\{\,\text{smashing subcategories of }\D(R)\,\big\}\\
\,\cong\, & \big\{\,\text{subsets of }\Spec\,R\text{ closed under specialization}\,\big\}
\end{align*}
of lattices.
\end{thm*}

Here, $\D(R)$ denotes the derived category of $\Mod\,R$.
A localizing subcategory is defined as a full triangulated subcategory which is closed under arbitrary direct sums, and a smashing subcategory is defined as a localizing subcategory such that Bousfield localization commutes with arbitrary direct sums.
The second main purpose of this paper is to construct a module version of the above Neeman's theorem:

\begin{thmb}
Let $R$ be a commutative noetherian ring.
Then one has an isomorphism
\begin{align*}
& \big\{\,\text{full subcategories of }\mod\,R\text{ closed under submodules and extensions}\,\big\}\\
\,\cong\, & \big\{\,\text{subsets of }\Spec\,R\,\big\}\\
\intertext{of lattices.
Moreover, this induces the isomorphism}
& \big\{\,\text{Serre subcategories of }\mod\,R\,\big\}\\
\,\cong\, & \big\{\,\text{subsets of }\Spec\,R\text{ closed under specialization}\,\big\}
\end{align*}
of lattices given in Theorem A.
\end{thmb}

In Section 2, we will give the precise definitions of subcategories which are stated above, and study several basic properties of them.
In Sections 3 and 4, we shall give proofs of Theorems A and B, respectively.

\begin{rem*}
It is known that Hovey's paper \cite{Hovey} contains an error, but it is not relevant to the results and arguments in this paper.
The error has recently been corrected by Krause \cite{Krause}.
\end{rem*}

\section{Basic properties}

In this section, we will give some definitions and several basic results most of which are necessary to state and prove the main results of this paper.

Let $\A$ be an additive category.
A full subcategory $\X$ of $\A$ is said to be {\em closed under isomorphisms} (or {\em replete}) provided that if $X$ is an object of $\X$ and $Y$ is an object of $\A$ which is isomorphic to $X$, then $Y$ is also an object of $\X$.
In this paper, by a {\em subcategory} we always mean a nonempty full subcategory which is closed under isomorphisms.

First of all, we recall the definitions of various types of closedness of a subcategory of an additive category.

\begin{defn}
Let $\A$ be an additive category, and let $\X$ be a full subcategory of $\A$.
We say that
\begin{enumerate}
\item[(1)]
$\X$ is {\em closed under subobjects} (resp. {\em closed under quotient objects}) provided that if $X$ is an object of $\X$ and $Y\in\A$ is a subobject (resp. a quotient object) of $X$, then $Y$ is also an object of $\X$.
\item[(2)]
$\X$ is {\em closed under direct summands} (or {\em closed under retracts}) provided that if $X$ is an object of $\X$ and $Y\in\A$ is a direct summand of $X$, then $Y$ is also an object of $\X$.
\item[(3)]
$\X$ is {\em closed under finite direct sums} (resp. {\em closed under arbitrary direct sums}) if all finite direct sums (resp. arbitrary direct sums) of objects of $\X$ are objects of $\X$.
\item[(4)]
$\X$ is {\em closed under extensions} provided that for any exact sequence $0 \to A \to B \to C \to 0$ in $\A$, if both $A$ and $C$ are objects of $\X$, then so is $B$.
\item[(5)]
$\X$ is {\em closed under kernels} (resp. {\em closed under images}, {\em closed under cokernels}) if the kernel (resp. the image, the cokernel) of every morphism of objects of $\X$ is also an object of $\X$.
\item[(6)]
$\X$ is {\em closed under homologies} if the homologies of every chain complex of objects of $\X$ are objects of $\X$.
\item[(7)]
$\X$ is {\em closed under direct limits} provided that if $\{ X_\lambda\}_{\lambda\in\Lambda}$ is a direct system of objects of $\X$ then the direct limit $\varinjlim _{\lambda\in\Lambda}X_\lambda$ is an object of $\X$.
\end{enumerate}
\end{defn}

Here we study the relationships among the closed properties of a subcategory which are defined above.

\begin{prop}\label{bbasic}
Let $\A$ be an abelian category and $\X$ a subcategory of $\A$.
Then the following hold.
\begin{enumerate}
\item[{\rm (1)}]
If $\X$ is closed under kernels or cokernels, then $\X$ is closed under direct summands and contains the zero object of $\A$.
\item[{\rm (2)}]
If $\X$ is closed under extensions, then $\X$ is closed under finite direct sums.
\item[{\rm (3)}]
If $\X$ is closed under kernels and cokernels, then $\X$ is closed under images and homologies.
\item[{\rm (4)}]
If $\X$ is closed under subobjects and cokernels, then $\X$ is closed under quotient objects.
\item[{\rm (5)}]
If $\X$ is closed under arbitrary direct sums and quotient objects, then $\X$ is closed under direct limits.
\end{enumerate}
\end{prop}

\begin{pf}
(1) Let $X$ be an object of $\X$ and $Y$ an object of $\A$ which is a direct summand of $X$.
Then we can write $X=Y\oplus Z$ for some object $Z$ of $\A$.
Considering the morphism $f: X \to X$ given by $(y,z) \mapsto (0,z)$ for $y\in Y$ and $z\in Z$, we see that both the kernel and the cokernel of $f$ are isomorphic to $Y$.
By the assumption that $\X$ is closed under kernels or cokernels, the object $Y$ is in $\X$.
Therefore $\X$ is closed under direct summands.

On the other hand, since $\X$ is nonempty, there exists an object $W\in\X$.
Let $g:W\to W$ be the identity morphism.
Then both the kernel and the cokernel of $g$ are the zero object.
Since $\X$ is closed under kernels or cokernels, $\X$ contains the zero object.

(2) Let $X$ and $Y$ be two objects of $\X$.
Then there exists a natural split exact sequnece
$$
0 \to X \to X\oplus Y \to Y \to 0
$$
in $\A$.
Since $\X$ is closed under extensions, we see from this exact sequence that the direct sum $X\oplus Y$ is an object of $\X$.
This argument shows that $\X$ is closed under finite direct sums.

(3) Let $f:X \to Y$ be a morphism of objects of $\X$.
Then there is an exact sequence
$$
0 \to \Im\,f \to Y \overset{\pi}{\to} \Coker\,f \to 0
$$
in $\A$.
Since $\X$ is closed under cokernels, the object $\Coker\,f$ is in $\X$.
Noting that the object $\Im\,f$ is the kernel of $\pi$ and $\X$ is closed under kernels, we see that $\Im\,f$ is in $\X$.
Therefore $\X$ is closed under images.
That $\X$ is closed under kernels and images implies that $\X$ is closed under homologies.

(4) Let $X$ be an object of $\X$ and $Y\in\A$ a quotient object of $X$.
Then there exists a subobject $Z\in\A$ of $X$ such that $Y=X/Z$.
Since $\X$ is closed under subobjects, $Z$ is an object of $\X$.
Let $i:Z\to X$ be the natural inclusion.
Since $Y$ coincides with the cokernel of $i$ and $\X$ is closed under cokernels, we have $Y\in\X$.
This says that $\X$ is closed under quotient objects.

(5) Let $\{ X_\lambda\}_{\lambda\in\Lambda}$ be a direct system of objects of $\X$, and let $X=\varinjlim _{\lambda\in\Lambda}X_\lambda$ be the direct limit.
Then, by definition, $X$ is a quotient object of $Y:=\bigoplus _{\lambda\in\Lambda}X_\lambda$.
Since $\X$ is closed under arbitrary direct sums, $Y$ is an object of $\X$.
Since $\X$ is closed under quotient objects, $X$ is an object of $\X$.
Consequently, $\X$ is closed under direct limits.
\qed
\end{pf}

Next, we recall the definitions of a coherent subcategory, a Serre subcategory and a torsion class, which will play important roles throughout this paper.

\begin{defn}
Let $\A$ be an abelian category, and let $\X$ be a subcategory of $\A$.
Then
\begin{enumerate}
\item[(1)]
$\X$ is called a {\em coherent subcategory} of $\A$ if it is closed under kernels, cokernels and extensions.
\item[(2)]
$\X$ is called a {\em Serre subcategory} of $\A$ if it is a coherent subcategory which is closed under subobjects.
\item[(3)]
$\X$ is called a {\em (hereditary) torsion class} of $\A$ if it is a Serre subcategory which is closed under arbitrary direct sums.
\end{enumerate}
\end{defn}

\begin{rem}
The original definition of a coherent subcategory is as follows: let $\A$ be an abelian category, and let $\X$ be a subcategory of $\A$.
It is said that $\X$ is coherent provided that for any exact sequence
$$
A \to B \to C \to D \to E
$$
in $\A$, if $A$, $B$, $D$ and $E$ are in $\A$, then so is $C$.
One can easily check that this definition is equivalent to our definition.
\end{rem}

A coherent subcategory, a Serre subcategory and a torsion class have the following properties, which immediately follow from Proposition \ref{bbasic}.

\begin{cor}\label{basic}
Let $\A$ be an abelian category.
\begin{enumerate}
\item[(1)]
Let $\X$ be a coherent subcategory of $\A$.
Then $\X$ contains the zero object of $\A$, and $\X$ is closed under finite direct sums, direct summands, images and homologies.
\item[(2)]
Let $\X$ be a Serre subcategory of $\A$.
Then $\X$ is closed under quotient objects.
\item[(3)]
Let $\X$ be a torsion class of $\A$.
Then $\X$ is closed under direct limits.
\end{enumerate}
\end{cor}

Throughout the rest of this paper, let $R$ be a commutative ring.
We denote by $\Mod\,R$ the category of $R$-modules, and by $\mod\,R$ the full subcategory of $\Mod\,R$ consisting of finitely presented $R$-modules.
Let us recall the definitions of a lattice and a homomorphism of lattices.

\begin{defn}
\begin{enumerate}
\item[(1)]
Let $\L$ be an ordered set.
\begin{enumerate}
\item[(i)]
Let $x,y\in\L$ be elements.
If the supremum (resp. infimum) of the set $\{ x,y\}$ exists, then it is called the {\em join} (resp. {\em meet}) of $x$ and $y$, and denoted by $x\vee y$ (resp. $x\wedge y$).
\item[(ii)]
It is said that $\L$ is a {\em lattice} if any two elements of $\L$ have both the join and the meet.
\end{enumerate}
\item[(2)]
A map $f:\L \to \L'$ of lattices is called a {\em (lattice) homomorphism} if $f(x\vee y)=f(x)\vee f(y)$ and $f(x\wedge y)=f(x)\wedge f(y)$ for all $x,y\in\L$.
A bijective lattice homomorphism is called a {\em (lattice) isomorphism}.
\end{enumerate}
\end{defn}

This paper will deal with the following lattices of subcategories of modules.

\begin{defn}
\begin{enumerate}
\item[(1)]
Let $\coh(R)$ be the coherent subcategory of $\Mod\,R$ generated by $R$.
We denote by $\lcoh(\coh(R))$ the lattice of all coherent subcategories of $\coh(R)$.
\item[(2)]
Let $\Serre(R)$ be the Serre subcategory of $\Mod\,R$ generated by $R$.
We denote by $\lserre(\Serre(R))$ the lattice of all Serre subcategories of $\Serre(R)$.
\item[(3)]
We denote by $\ltors(\Mod\,R)$ the lattice of all torsion classes of $\Mod\,R$.
\end{enumerate}
\end{defn}

\begin{rem}\label{infact}
If $R$ is noetherian, then one has $\coh(R) = \Serre(R) = \mod\,R$; see the first sentence and the latter half of Page 3185 of \cite{Hovey}.
Hence, whenever $R$ is a noetherian ring, one has $\lcoh(\coh(R))=\lcoh(\mod\,R)$ and $\lserre(\Serre(R))=\lserre(\mod\,R)$.
\end{rem}

A {\it perfect} $R$-complex $P_\bullet$ is defined to be an $R$-complex of the form
$$
P_\bullet=(0 \to P_s \to P_{s-1} \to \cdots \to P_{t+1} \to P_t \to 0), 
$$
where each $P_i$ is a finitely generated projective $R$-module.
We denote by $\D (R)$ the derived category of the category $\Mod\,R$, and by $\dperf(R)$ the full subcategory of $\D(R)$ consisting of $R$-complexes which are isomorphic to perfect $R$-complexes.

\begin{rem}
Recall that an object $X$ of an additive category $\A$ is called {\em small} (or {\em compact}) if the functor $\Hom _\A(X,-)$ preserves arbitrary direct sums.
It is well-known that a complex of $R$-modules is quasi-isomorphic to a perfect complex if and only if it is a small object of $\D(R)$ (cf. \cite[3.7]{DGI}).
\end{rem}

\begin{defn}
\begin{enumerate}
\item[(1)]
Let $\T$ be a triangulated category and $\X$ a subcategory of $\T$.
Then we say that $\X$ is a {\it thick subcategory} if it satisfies the following two conditions.
\begin{enumerate}
\item[(a)]
$\X$ is closed under direct summands.
\item[(b)]
For any exact triangle $A \to B \to C \to \Sigma A$ in $\T$, if two of the objects $A,B,C$ are in $\X$, then so is the third.
\end{enumerate}
\item[(2)]
We denote by $\lthick (\dperf(R))$ the lattice of all of the thick subcategories of $\D (R)$ whose objects are small, namely, the lattice of all thick subcategories of $\dperf(R)$.
\end{enumerate}
\end{defn}

Recall that a subset $S$ of $\Spec\,R$ is said to be {\em closed under specialization} provided that if $\p$ is a prime ideal in $S$ and $\q$ is a prime ideal containing $\p$ then $\q$ is also in $S$.
Dually, $S$ is said to be {\em closed under generalization} provided that if $\p$ is a prime ideal in $S$ and $\q$ is a prime ideal contained in $\p$ then $\q$ is also in $S$.
Note that every union of closed subsets of $\Spec\,R$ is closed under specialization.
Similarly, every intersection of open subsets of $\Spec\,R$ is closed under generalization.

\begin{defn}
We denote by $\lsp(\Spec\,R)$ the lattice of all subsets of $\Spec\,R$ that are closed under specialization, and by $\lsp^0(\Spec\,R)$ the sublattice of $\lsp(\Spec\,R)$ consisting of all complements of arbitrary intersections of quasi-compact open subsets of $\Spec\,R$.
\end{defn}

\begin{rem}\label{compact}
An open subset $U$ of $\Spec\,R$ is quasi-compact if and only if $U=D(I):=\{\p\in\Spec\,R\mid I\nsubseteq\p\}$ for some finitely generated ideal $I$ of $R$; see the argument on the top of \cite[Page 72]{Hartshorne}.
Therefore, if $R$ is noetherian, then every open subset of $\Spec\,R$ is quasi-compact, and $\lsp^0(\Spec\,R)$ consists of all unions of closed subsets of $\Spec\,R$.
\end{rem}

Hovey \cite{Hovey} constructs the following order-preserving maps among the lattices which we defined above:
$$
\begin{CD}
\lsp(\Spec\,R) \\
\begin{smallmatrix}
@A{\tau}AA & @VV{\sigma}V
\end{smallmatrix}\\
\ltors(\Mod\,R) \\
\begin{smallmatrix}
@A{\nu}AA & @VV{\mu}V
\end{smallmatrix}\\
\lserre(\Serre(R)) \\
\begin{smallmatrix}
@A{\beta}AA & @VV{\alpha}V
\end{smallmatrix}\\
\lcoh(\coh(R)) \\
\begin{smallmatrix}
@A{g}AA & @VV{f}V
\end{smallmatrix}\\
\lthick(\dperf(R))
\end{CD}
$$
The above maps are defined as follows.
\begin{enumerate}
\item[(1)]
For $S\in \lsp(\Spec\,R)$, let $\sigma(S)$ be the full subcategory of $\Mod\,R$ consisting of all $R$-modules $M$ with $\Supp\,M\subseteq S$.
For $\X\in\ltors (\Mod\,R)$, let $\tau(\X)$ be the union $\bigcup_{M\in\X}\Supp\,M$.
\item[(2)]
For $\X\in\ltors (\Mod\,R)$, let $\mu(\X)$ be the full subcategory of $\Mod\,R$ consisting of all $R$-modules $M$ in the intersection of $\X$ and $\Serre(R)$.
For $\Y\in\lserre (\Serre(R))$, let $\nu(\Y)$ be the torsion class of $\Mod\,R$ generated by $\Y$.
\item[(3)]
For $\Y\in\lserre (\Serre(R))$, let $\alpha (\Y)$ be the full subcategory of $\Mod\,R$ consisting of all $R$-modules $M$ in the intersection of $\Y$ and $\coh(R)$.
For $\zzz\in\lcoh (\coh(R))$, let $\beta (\zzz)$ be the Serre subcategory of $\Mod\,R$ generated by $\zzz$.
\item[(4)]
For $\zzz\in\lcoh (\coh(R))$, let $f(\zzz)$ be the full subcategory of $\D (R)$ consisting of all complexes $X_\bullet\in\dperf(R)$ such that the $i$th homology $H_i(X_\bullet)$ is in $\zzz$ for any $i\in\Z$.
For $\W\in\lthick (\dperf(R))$, let $g(\W)$ be the coherent subcategory of $\Mod\,R$ generated by $\{\, H_i(X_\bullet)\,|\,X_\bullet\in\W,i\in\Z\,\}$.
\end{enumerate}

We recall here the definition of an adjoint pair of order-preserving maps.
Let $\phi: A \to B$ and $\psi: B \to A$ be two order-preserving maps between ordered sets.
Then the pair $(\phi,\psi)$ is said to be an {\em adjoint pair} provided that $\phi(a)\le b$ if and only if $a\le \psi(b)$ for any $a\in A$ and $b\in B$.
Concerning the above order-preserving maps among lattices, the following proposition holds.

\begin{prop}\label{morphisms}
\begin{enumerate}
\item[{\rm (1)}]
\begin{enumerate}
\item[{\rm (a)}]
The pair $(\tau,\sigma)$ is an adjoint pair.
\item[{\rm (b)}]
The composite map $\tau\sigma$ is the identity map.
\item[{\rm (c)}]
If $R$ is noetherian, then $\sigma$ is a lattice isomorphism and $\tau$ is the inverse homomorphism.
\end{enumerate}
\item[{\rm (2)}]
\begin{enumerate}
\item[{\rm (a)}]
The pair $(\nu,\mu)$ is an adjoint pair.
\item[{\rm (b)}]
The composite map $\nu\mu$ is the identity map.
\item[{\rm (c)}]
If $R$ is noetherian, then $\mu$ is a lattice isomorphism and $\nu$ is the inverse homomorphism.
\end{enumerate}
\item[{\rm (3)}]
\begin{enumerate}
\item[{\rm (a)}]
The pair $(\beta,\alpha)$ is an adjoint pair.
\item[{\rm (b)}]
If $R$ is noetherian, then the composite map $\beta\alpha$ is the identity map.
\end{enumerate}
\item[{\rm (4)}]
\begin{enumerate}
\item[{\rm (a)}]
The pair $(g,f)$ is an adjoint pair.
\item[{\rm (b)}]
The composite map $fg$ is the identity map.
\end{enumerate}
\end{enumerate}
\end{prop}

\begin{pf}
It is easy to check that every order-preserving bijective map between two lattices is a lattice isomorphism.
Hence we see from the arguments in Page 3185 of \cite{Hovey} that the assertions (1), (2) and (3) hold.
(Here, note that derived categories do not appear in those arguments.)

The assertion (a) and (b) in (4) are shown in \cite[Proposition 1.4]{Hovey} and \cite[Corollary 2.2]{Hovey}, respectively.
\qed
\end{pf}

The above proposition especially says that one has the following relationships between two modules whose supports have inclusion relation.

\begin{cor}\label{generated}
Let $R$ be a noetherian ring, and let $M,N$ be $R$-modules with $\Supp\,M\subseteq\Supp\,N$.
Then $M$ belongs to the torsion class of $\Mod\,R$ generated by $N$.
If $M$ and $N$ are finitely generated, then $M$ belongs to the Serre subcategory of $\mod\,R$ generated by $N$.
\end{cor}

\begin{pf}
Proposition \ref{morphisms} says that the maps $\sigma,\mu$ are isomorphisms whose inverse maps are $\tau,\nu$ respectively.
Let $\T$ be the torsion class of $\Mod\,R$ generated by $N$.
We have $\T=\sigma\tau(\T)$, which is the full subcategory of $\Mod\,R$ consisting of all $R$-modules $K$ with $\Supp\,K\subseteq\bigcup_{L\in\T}\Supp\,L$.
Since $N$ is in $\T$, we get $\Supp\,M\subseteq\Supp\,N\subseteq\bigcup_{L\in\T}\Supp\,L$.
Therefore $M$ is in $\T$.

Suppose that both $M$ and $N$ are finitely generated.
Let $\sss$ be the Serre subcategory of $\mod\,R$ generated by $N$.
Then we have $\sss=\mu\sigma\tau\nu(\sss)$, which consists of all finitely generated $R$-modules $K$ with $\Supp\,K\subseteq\bigcup_{L\in\T}\Supp\,L(=\bigcup_{L\in\sss}\Supp\,L)$.
Hence $M$ is in $\sss$.
\qed
\end{pf}

\section{Coherent subcategories are Serre}

Throughout this section, let $R$ be a commutative ring.
We begin with proving the following theorem.

\begin{thm}\label{main1}
Let $R$ be a noetherian ring.
Let $\X$ be a full subcategory of $\mod\,R$ which is closed under finite direct sums, kernels and cokernels.
Then $\X$ is closed under submodules and quotient modules.
\end{thm}

\begin{pf}
According to Proposition \ref{bbasic}(4), it is enough to prove that $\X$ is closed under submodules.
Assume that $\X$ is not closed under submodules.
Then there exist an $R$-module $X$ in $\X$ and an $R$-submodule $M$ of $X$ such that $M$ does not belong to $\X$.
Since $R$ is noetherian and $X$ is a finitely generated $R$-module, $X$ is a noetherian $R$-module.
Hence we can choose $M$ to be a maximal element, with respect to the inclusion relation, of the set of $R$-submodules $M'$ of $X$ such that $M'$ does not belong to $\X$.
Since $M$ does not coincide with $X$, there is an element $x\in X-M$.
Set $Y=M+Rx$.
Note that $Y$ is an $R$-submodule of $X$ strictly containing $M$.
By the maximality of $M$, the module $Y$ is in $\X$.
Put $I=(M:x):=\{ a\in R\mid ax\in M \}$.
This is an ideal of $R$, and we easily see that the quotient $R$-module $Y/M$ is isomorphic to $R/I$.
There is an exact sequence
$$
0 \to M \to Y \overset{\pi}{\to} R/I \to 0
$$
of $R$-modules.
Since $M\notin\X$ and $Y\in\X$ and $\X$ is closed under kernels, we see from this exact sequence that $R/I$ must not be in $\X$.

On the other hand, the map $\pi$ in the exact sequence induces a surjective homomorphism
$$
\overline\pi :Y/IY \to R/I
$$
of $R/I$-modules, which sends the residue class of $y\in Y$ in $Y/IY$ to $\pi (y)$.
Of course $R/I$ is a projective $R/I$-module, so $\overline\pi$ is a split epimorphism.
Therefore $R/I$ is isomorphic to a direct summand of $Y/IY$.
The noetherian property of $R$ implies that the ideal $I$ is finitely generated; write $I=(a_1,a_2,\dots,a_n)R$ for some elements $a_1,a_2,\dots,a_n\in R$.
There is an exact sequence
$$
\begin{CD}
R^{\oplus n} @>(a_1,\dots,a_n)>> R @>>> R/I @>>> 0
\end{CD}
$$
of $R$-modules.
Tensoring the $R$-module $Y$ with this exact sequence yields another exact sequence of $R$-modules:
$$
\begin{CD}
Y^{\oplus n} @>(a_1,\dots,a_n)>> Y @>>> Y/IY @>>> 0.
\end{CD}
$$
The assumption of the theorem says that $\X$ is closed under finite direct sums, cokernels, and direct summands; see Proposition \ref{bbasic}(1).
Hence the direct sum $Y^{\oplus n}$ belongs to $\X$, and so does the module $Y/IY$, and therfore so does $R/I$.
This is a contradiction, which says that $\X$ is closed under submodules.
Thus the proof of the theorem is completed.
\qed
\end{pf}

Corollary \ref{basic}(1) says that any coherent subcategory of $\mod\,R$ is closed under finite direct sums, kernels and cokernels.
Thus, according to Theorem \ref{main1}, we obtain the following result, which is the former half part of Theorem A in the first section of this paper.

\begin{cor}\label{wideserre}
Let $R$ be a noetherian ring.
Then every coherent subcategory of $\mod\,R$ is a Serre subcategory of $\mod\,R$.
\end{cor}

Now, let us check that the subset $\lsp(\Spec\,R)$ coincides with $\lsp^0(\Spec\,R)$ if $R$ is a noetherian ring.

\begin{prop}\label{J}
\begin{enumerate}
\item[{\rm (1)}]
Let $Z$ be a subset of $\Spec\,R$ which is closed under specialization.
Then one has
$$
Z=\bigcup_{\p\in Z}V(\p).
$$
\item[{\rm (2)}]
Let $R$ be a noetherian ring.
Then
$$
\lsp(\Spec\,R)=\lsp^0(\Spec\,R).
$$
\end{enumerate}
\end{prop}

\begin{pf}
(1) Let $\q$ be a prime ideal in $Z$.
Then $\q$ is in $V(\q)$, which is contained in $\bigcup_{\p\in Z}V(\p)$.
As to the opposite inclusion relation, take a prime ideal $\q$ in $\bigcup_{\p\in Z}V(\p)$.
Then $\q$ is in $V(\p)$ for some $\p\in Z$.
Since $Z$ is closed under specialization, we get $\q\in Z$, as required.

(2) This immediately follows from Remark \ref{compact} and the assertion (1).
\qed
\end{pf}

To prove our next result, we prepare here the following two lemmas.

\begin{lem}\label{gf}
Let $R$ be a noetherian ring and $\M$ a subcategory of $\mod\,R$ which is closed under finite direct sums and cokernels.
Let $M$ be a cyclic $R$-module in $\M$.
Then there exists a perfect $R$-complex $X_\bullet$ such that $H_0(X_\bullet)$ is isomorphic to $M$ and that $H_j(X_\bullet)$ belongs to $\M$ for any $j\in\Z$.
\end{lem}

\begin{pf}
Since $M$ is cyclic, there exists an ideal $I$ of $R$ such that $M$ is isomorphic to $R/I$.
The noetherian property of $R$ implies that the ideal $I$ is finitely generated; let $\xx = x_1,x_2,\dots,x_r$ be a system of generators of $I$.
Consider the Koszul complex $K_\bullet := K_\bullet (\xx , R)$ of the sequence $\xx$.
The complex $K_\bullet$ is a perfect $R$-complex, and the zeroth homology $H_0(K_\bullet)$ is equal to $R/I$, which is isomorphic to $M$.
Thus, to show the lemma, it suffices to check that the homology $H_j(K_\bullet)$ belongs to $\M$ for each $j\in\Z$.
Note from \cite[Proposition 1.6.5(b)]{BH} that the $R$-module $H_j(K_\bullet)=H_j(\xx,R)$ is annihilated by the ideal $I=\xx R$.
Hence $H_j(K_\bullet)$ can be regarded as an $R/I$-module.
Since $R$ is noetherian, $H_j(K_\bullet)$ is finitely generated, hence finitely presented as an $R/I$-module.
It follows that there is an exact sequence
$$
(R/I)^{\oplus n} \to (R/I)^{\oplus m} \to H_j(K_\bullet) \to 0
$$
of $R/I$-modules.
Since $\M$ is closed under finite direct sums, the sum $(R/I)^{\oplus i}\cong M^{\oplus i}$ is an object of $\M$ for any $i\ge 0$.
Since $\M$ is closed under cokernels, the module $H_j(K_\bullet)$ is in $\M$, as desired.
\qed
\end{pf}

\begin{lem}\label{filt}
Let $M$ be a finitely generated $R$-module.
Then there exist exact sequences
$$
\begin{cases}
0 \to R/I_1 \to M_0 \to M_1 \to 0,\\
0 \to R/I_2 \to M_1 \to M_2 \to 0,\\
\qquad \vdots \\
0 \to R/I_{n-1} \to M_{n-2} \to M_{n-1} \to 0,\\
0 \to R/I_n \to M_{n-1} \to M_n \to 0
\end{cases}
$$
of $R$-modules such that $I_1,I_2, \dots , I_{n-1}, I_n$ are ideals of $R$, and $M_0=M$ and $M_n=0$.
\end{lem}

\begin{pf}
Let $x_1,x_2,\dots,x_n$ be a system of generators of $M$.
Set $M_0=M$.
We have an isomorphism $Rx_1\cong R/I_1$, where $I_1=\Ann(x_1)$.
Putting $M_1=M_0/Rx_1$, we have an exact sequence
$$
0 \to R/I_1 \to M_0 \to M_1 \to 0
$$
of $R$-modules.
Note that the $R$-module $M_1$ is generated by $n-1$ elements.
By induction on $n$, we can obtain such a system of exact sequences as in the lemma.
\qed
\end{pf}

Now we are in a position to prove the following theorem, which is the latter half part of Theorem A in the first section.
In the proof, we should note that all the isomorphisms in the theorem except $\lcoh(\coh(R))\cong\lthick(\dperf(R))$ are obtained without using derived categories.

\begin{thm}
Let $R$ be a noetherian ring.
Then the homomorphisms $\sigma,\mu,\alpha,f$ (defined in the previous section) are lattice isomorphisms, and $\tau,\nu,\beta,g$ are their inverse homomorphisms, respectively.
Consequently, one has
\begin{align*}
\lsp^0(\Spec\,R) & = \lsp(\Spec\,R) \\
& \cong \ltors(\Mod\,R) \\
& \cong \lserre(\Serre(R)) \\
& = \lcoh(\coh(R)) \\
& \cong \lthick(\dperf(R)).
\end{align*}
\end{thm}

\begin{pf}
The equality $\lsp^0(\Spec\,R) = \lsp(\Spec\,R)$ is already shown in Proposition \ref{J}.
Proposition \ref{morphisms} says that the homomorphisms $\sigma,\mu$ are lattice isomorphisms and $\tau,\nu$ are the inverse homomorphisms of $\sigma,\mu$ respectively.
Hence we have $\lsp(\Spec\,R)\cong \ltors(\Mod\,R)\cong \lserre(\Serre(R))$.
It is seen from Corollary \ref{wideserre} and Remark \ref{infact} that both of the homomorphisms $\alpha$ and $\beta$ are the identity maps, and we have $\lserre(\Serre(R))= \lcoh(\coh(R))$.

It remains to prove that $f$ is an isomorphism with the inverse homomorphism $g$.
The composite map $fg$ is the identity homomorphism by Proposition \ref{morphisms}(4)(b), and the subcategory $gf(\X)$ is contained in $\X$ for every $\X\in\lcoh (\coh(R))$ by Proposition \ref{morphisms}(4)(a).
Let us show that $\X$ is contained in $gf(\X)$.
Let $M$ be an $R$-module in $\X$.
Since Remark \ref{infact} guarantees that $\X$ is a subcategory of $\mod\,R$, $M$ is a finitely generated $R$-module, and we have a system of exact sequences as in Lemma \ref{filt}.
(In the following, we use the same notation.)
Since $\X$ is a coherent subcategory of $\mod\,R$ by Remark \ref{infact}, we see from Corollary \ref{wideserre} that $\X$ is a Serre subcategory of $\mod\,R$.
Hence $\X$ is closed under submodules and quotient modules in $\mod\,R$ by Corollary \ref{basic}(2).
From the above exact sequences we easily see that the cyclic $R$-module $R/I_i$ belongs to $\X$ for every $1\le i\le n$.
Note by Corollary \ref{basic}(1) that $\X$ is closed under finite direct sums.
Hence Lemma \ref{gf} shows that for each integer $1\le i\le n$ there exists a perfect $R$-complex $X_\bullet^{(i)}$ such that $R/I_i$ is isomorphic to $H_0(X_\bullet^{(i)})$ and that $H_j(X_\bullet^{(i)})$ belongs to $\X$ for any $j\in\Z$.
It easily follows from the definitions of the homomorphisms $f,g$ that the $R$-module $R/I_i$ belongs to $gf(\X)$.
Since $gf(\X)$ is a coherent subcategory, it is closed under extensions.
Hence from the system of exact sequences, we see that $M$ belongs to $gf(\X)$.
Therefore $\X $ is contained in $gf(\X )$, and thus $gf$ is the identity homomorphism.
This completes the proof of our theorem.
\qed
\end{pf}

\section{In relation to Neeman's classification}

In this section, we shall give a module version of Neeman's classification theorem of localizing categories and smashing categories, which we stated in the first section of this paper.
Throughout this section, let $R$ be a commutative noetherian ring.

We denote by $\lse(\mod\,R)$ the lattice of all subcategories of $\mod\,R$ which are closed under submodules and extensions, and by $\L(\Spec\,R)$ the lattice of all subsets of $\Spec\,R$.
We define maps
$$
\begin{cases}
\Phi:\L(\Spec\,R) \to \lse(\mod\,R) \\
\Psi:\lse(\mod\,R) \to \L(\Spec\,R)
\end{cases}
$$
by $\Phi(S)=\{\,M\in\mod\,R\mid\Ass\,M\subseteq S\,\}$ and $\Psi(\M)=\bigcup_{M\in\M}\Ass\,M$.
It is easy to check that these maps are lattice homomorphisms.
Let $\phi:\lsp(\Spec\,R) \to \lserre(\mod\,R)$ and $\psi:\lserre(\mod\,R) \to \lsp(\Spec\,R)$ be the composite maps $\mu\sigma$ and $\tau\nu$, where $\sigma$, $\tau$, $\mu$, $\nu$ are the homomorphisms defined in Section 2.
Note that the maps $\phi,\psi$ are given by $\phi(S)=\{\,M\in\mod\,R\mid\Supp\,M\subseteq S\,\}$ and $\psi(\M)=\bigcup_{M\in\M}\Supp\,M$ for $S\in\lsp(\Spec\,R)$ and $\M\in\lserre(\mod\,R)$.

Recall that $\phi$ is an isomorphism and $\psi$ is its inverse homomorphism since $R$ is noetherian.
This section is mainly devoted to proving the following theorem.

\begin{thm}\label{main}
Let $R$ be a commutative noetherian ring.
Then the homomorphism $\Phi$ is an isomorphism and $\Psi$ is its inverse homomorphism.
Moreover, $\Phi$ and $\Psi$ induce the isomorphisms $\phi$ and $\psi$, respectively.
Thus one has the following commutative diagram.
$$
\begin{CD}
\lse(\mod\,R) @>{\Psi}>{\cong}> \L(\Spec\,R) @>{\Phi}>{\cong}> \lse(\mod\,R) \\
@AA{\subseteq}A @AA{\subseteq}A @AA{\subseteq}A \\
\lserre(\mod\,R) @>{\psi}>{\cong}> \lsp(\Spec\,R) @>{\phi}>{\cong}> \lserre(\mod\,R)
\end{CD}
$$
\end{thm}

This theorem will be proved after showing the following two lemmas.

\begin{lem}\label{coprimary}
Let $\M$ be a subcategory of $\mod\,R$ which is closed under submodules and extensions, and let $M$ be a finitely generated $R$-module.
Suppose that $M$ has a unique associated prime $\p$.
If $R/\p$ is in $\M$, then so is $M$.
\end{lem}

\begin{pf}
Assume that $M$ is not in $\M$.
Set $M_0=M$, and let $f_{0,1},\dots,f_{0,s_0}$ be a system of generators of the $R$-module $\Hom_R(M_0,R/\p)$.
There is an exact sequence
$$
\begin{CD}
0 @>>> M_1 @>>> M_0 @>{\left(
\begin{smallmatrix}
f_{0,1}\\
\vdots\\
\\
f_{0,s_0}
\end{smallmatrix}
\right)}>> (R/\p)^{\oplus s_0}
\end{CD}
$$
of $R$-modules.
Since $M_0=M$ is not in $\M$ and $(R/\p)^{\oplus s_0}$ is in $\M$ and $\M$ is closed under submodules and extensions, it is easily seen that $M_1$ must not be in $\M$.
In particular, $M_1\ne 0$ and hence $\p$ is the unique associated prime of $M_1$.
Letting $f_{1,1},\dots,f_{1,s_1}$ be a system of generators of the $R$-module $\Hom_R(M_1,R/\p)$, we have an exact sequence
$$
\begin{CD}
0 @>>> M_2 @>>> M_1 @>{\left(
\begin{smallmatrix}
f_{1,1}\\
\vdots\\
\\
f_{1,s_1}
\end{smallmatrix}
\right)}>> (R/\p)^{\oplus s_1}.
\end{CD}
$$
Since $M_1$ is not in $\M$ and $(R/\p)^{\oplus s_1}$ is in $\M$, we see that $M_2$ is not in $\M$, and that $\p$ is the unique associated prime of $M_2$.
Iterating this procedure, for each integer $i\ge 0$ we obtain an exact sequence
$$
\begin{CD}
0 @>>> M_{i+1} @>>> M_i @>{\left(
\begin{smallmatrix}
f_{i,1}\\
\vdots\\
\\
f_{i,s_i}
\end{smallmatrix}
\right)}>> (R/\p)^{\oplus s_i},
\end{CD}
$$
where $f_{i,1},\dots,f_{i,s_i}$ is a system of generators of the $R$-module $\Hom_R(M_i,R/\p)$ and $\p$ is the unique associated prime of $M_i$.
Localizing the descending chain $M=M_0\supseteq M_1\supseteq M_2\supseteq\cdots$ at $\p$ yields a descending chain
$$
M_\p=(M_0)_\p\supseteq (M_1)_\p\supseteq (M_2)_\p\supseteq\cdots
$$
of $R_\p$-modules.
Since the $R_\p$-module $(M_i)_\p$ has finite length for every $i$, there exists an integer $t$ such that $(M_t)_\p=(M_{t+1})_\p=(M_{t+2})_\p=\cdots$.
The exact sequence
$$
\begin{CD}
0 @>>> (M_{t+1})_\p @>{=}>> (M_t)_\p @>{\left(
\begin{smallmatrix}
(f_{t,1})_\p\\
\vdots\\
\\
(f_{t,s_t})_\p
\end{smallmatrix}
\right)}>> \kappa(\p)^{\oplus s_t},
\end{CD}
$$
shows that $\Hom_{R_\p}((M_t)_\p,\kappa(\p))=R_\p(f_{t,1})_\p +\cdots +R_\p(f_{t,s_t})_\p =0$.
Therefore $(M_t)_\p =0$.
This is a contradiction since $\p\in\Ass\,M_t\subseteq\Supp\,M_t$.
Thus we conclude that $M$ is in $\M$.
\qed
\end{pf}

\begin{lem}\label{seisitu}
Let $\M$ be a subcategory of $\mod\,R$ which is closed under submodules and extensions.
Let $M$ be a finitely generated $R$-module.
Suppose that $R/\p$ belongs to $\M$ for every $\p\in\Ass\,M$.
Then $M$ also belongs to $\M$.
\end{lem}

\begin{pf}
Let $\p_1,\dots,\p_s$ be the associated primes of $M$, and let
$$
0 = N_1\cap\cdots\cap N_s
$$
be an irredundant primary decomposition of the zero submodule $0$ of $M$, where $N_i$ is a $\p_i$-primary submodule of $M$ for $1\le i\le s$.
Then the natural homomorphism
$$
M=M/N_1\cap\cdots\cap N_s\to M/N_1\oplus\cdots\oplus M/N_s
$$
is injective.
Since $\p_i$ is the unique associated prime of the $R$-module $M/N_i$, Lemma \ref{coprimary} implies that $M/N_i$ belongs to $\M$ for $1\le i\le s$.
Hence $M/N_1\oplus\cdots\oplus M/N_s$ belongs to $\M$, and so does $M$.
\qed
\end{pf}

Now we can achieve the main purpose of this section.

\begin{tpf}
Let $S$ be a subset of $\Spec\,R$.
The set $\Psi\Phi(S)$ is the union of $\Ass\,M$ where $M$ runs through finitely generated $R$-modules all of whose associated primes are in $S$.
It is trivial that this set is contained in $S$.
For a prime ideal $\p\in S$, we have $\Ass_R\,R/\p =\{\p\}\subseteq S$.
Hence $\p$ belongs to $\Psi\Phi(S)$, and therefore $\Psi\Phi(S)=S$.
Let $\M$ be a subcategory of $\mod\,R$ which is closed under submodules and extensions.
We have that $\Phi\Psi(\M)$ is the subcategory of $\mod\,R$ consisting of all finitely generated $R$-modules $N$ with $\Ass\,N\subseteq\bigcup_{M\in\M}\Ass\,M$, and it is obvious that $\Phi\Psi(\M)$ contains $\M$.
Let $N$ be a finitely generated $R$-module with $\Ass\,N\subseteq\bigcup_{M\in\M}\Ass\,M$.
Fix a prime ideal $\p\in\Ass\,N$.
Then there exists an $R$-module $M\in\M$ with $\p\in\Ass\,M$.
There is an injective homomorphism $R/\p\to M$, and $R/\p$ belongs to $\M$ since $\M$ is closed under submodules.
It follows from Lemma \ref{seisitu} that $N$ is in $\M$.
Hence $\Phi\Psi(\M)=\M$.
Thus we conclude that $\Phi$ is an isomorphism whose inverse homomorphism is $\Psi$.

On the other hand, let $S$ be a subset of $\Spec\,R$ which is closed under specialization.
Let $M$ be a finitely generated $R$-module such that $\Ass\,M$ is contained in $S$, and take $\p\in\Supp\,M$.
Then there is a prime ideal $\q\in\Min\,M\subseteq\Ass\,M$ that is contained in $\p$.
Since $\q$ is in $S$ and $S$ is closed under specialization, $\p$ is also in $S$.
Thus $\Phi(S)=\{\,M\in\mod\,R\mid\Ass\,M\subseteq S\,\}$ coincides with $\phi(S)=\{\,M\in\mod\,R\mid\Supp\,M\subseteq S\,\}$.
Let $\M$ be a Serre subcategory of $\mod\,R$.
Let $N\in\M$ and $\p\in\Supp\,N$.
Choose a prime ideal $\q\in\Min\,N$ which is contained in $\p$.
Then $\q$ is an associated prime of $N$, so there is an injective homomorphism $R/\q\to N$.
Since $\M$ is closed under submodules, the module $R/\q$ is in $\M$.
Noting that there is a surjective homomorphism $R/\q \to R/\p$ and that $\M$ is closed under quotient modules, $R/\p$ is also in $\M$.
Hence we get $\p\in\Ass\,R/\p\subseteq\bigcup_{M\in\M}\Ass\,M$.
Therefore the set $\psi(\M)=\bigcup_{M\in\M}\Supp\,M$ is contained in $\Psi(\M)=\bigcup_{M\in\M}\Ass\,M$, and we see that $\Psi(\M)=\psi(\M)$.
It follows that $\phi$ and $\psi$ are induced from $\Phi$ and $\Psi$, respectively.
\qed
\end{tpf}

Here, let us check that an analogous result to Corollary \ref{generated} holds.
This actually follows from Lemma \ref{seisitu}.

\begin{cor}
Let $R$ be a noetherian ring.
Let $M$ and $N$ be finitely generated $R$-modules with $\Ass\,M\subseteq\Ass\,N$.
Then $M$ is in the full subcategory of $\mod\,R$ closed under submodules and extensions which is generated by $N$.
\end{cor}

\begin{pf}
Let $\E$ be the full subcategory of $\mod\,R$ closed under submodules and extensions which is generated by $N$.
According to Lemma \ref{seisitu}, we have only to show that the $R$-module $R/\p$ is in $\E$ for every $\p\in\Ass\,M$.
Let $\p$ be a prime ideal in $\Ass\,M$.
The assumption says that $\p$ is in $\Ass\,N$.
Hence there exists an injective homomorphism $R/\p\to N$ of $R$-modules.
Since $N$ is in $\E$ and $\E$ is closed under submodules, $R/\p$ is also in $\E$, as required.
\qed
\end{pf}

In the following example, we will give several correspondences between subcategories of $\mod\,R$ which are closed under submodules and extensions and subsets of $\Spec\,R$, which are made by the isomorphisms $\Phi$ and $\Psi$.
Before that, we need to prepare some notation.
Let $I$ be an ideal of $R$, and let $M,N$ be $R$-modules.
We denote by $\Gamma_I(M)$ the {\it $I$-torsion submodule} of $M$, namely, the set of elements of $M$ which are annihilated by some power of $I$.
Recall that an $R$-module $M$ is called {\it $I$-torsion} if $\Gamma_I(M)=M$, and that $M$ is called {\it $I$-torsionfree} if $\Gamma_I(M)=0$.
It is well-known and easy to see that $M$ is $I$-torsion if and only if $\Ass\,M\subseteq V(I)$, and that $M$ is $I$-torsionfree if and only if $\Ass\,M\cap V(I)=\emptyset$.
We set $\grade(N,M)=\inf\{\,i\,|\,\Ext_R^i(N,M)\ne 0\,\}$, $\grade(I,M)=\grade(R/I,M)$, $\grade\,I=\grade(I,R)$ and $\grade\,M=\grade(\Ann\,M,R)$.

\begin{ex}
The isomorphisms $\Phi$ and $\Psi$ make the following correspondences.
Let $n$ be a nonnegative integer, $I$ an ideal of $R$ and $X$ a finitely generated $R$-module.
\begin{enumerate}
\item[{\rm (1)}]
$\big\{\,M\in\mod\,R\ \big|\ M\text{ is }I\text{-torsion}\,\big\}\,\leftrightarrow\,V(I)$.\\
\item[{\rm (2)}]
$\big\{\,M\in\mod\,R\ \big|\ \grade(X,M)>0\,\big\}\,\leftrightarrow\,\Spec\,R\setminus\Supp\,X$.\\
\item[{\rm (3)}]
$\big\{\,M\in\mod\,R\ \big|\ M\text{ is }I\text{-torsionfree}\,\big\}\\
\,=\,\big\{\,M\in\mod\,R\ \big|\ \grade(I,M)>0\,\big\}\,\leftrightarrow\,D(I)$.\\
\item[{\rm (4)}]
$\big\{\,M\in\mod\,R\ \big|\ \grade(M,X)\ge n\,\big\}\,\leftrightarrow\,\big\{\,\p\in\Spec\,R\ \big|\ \grade(\p,X)\ge n\,\big\}$.\\
\item[{\rm (5)}]
$\big\{\,M\in\mod\,R\ \big|\ \rank\,M=0\,\big\}\,=\,\big\{\,M\in\mod\,R\ \big|\ \grade\,M>0\,\big\}\\
\,\leftrightarrow\,\big\{\,\p\in\Spec\,R\ \big|\ \grade\,\p>0\,\big\}$.\\
\item[{\rm (6)}]
$\big\{\,M\in\mod\,R\ \big|\ \text{every }X\text{-regular element is }M\text{-regular}\,\big\}\\
\,\leftrightarrow\,\big\{\,\p\in\Spec\,R\ \big|\ \grade(\p,X)=0\,\big\}$.\\
\item[{\rm (7)}]
$\big\{\,M\in\mod\,R\ \big|\ M\text{ is torsionfree}\,\big\}\,\leftrightarrow\,\big\{\,\p\in\Spec\,R\ \big|\ \grade\,\p=0\,\big\}$.\\
\item[{\rm (8)}]
$\big\{\,M\in\mod\,R\ \big|\ \height\,\Ann\,M\ge n\,\big\}\,\leftrightarrow\,\big\{\,\p\in\Spec\,R\ \big|\ \height\,\p\ge n\,\big\}$.\\
\item[{\rm (9)}]
$\big\{\,M\in\mod\,R\ \big|\ \kdim\,M\le n\,\big\}\,\leftrightarrow\,\big\{\,\p\in\Spec\,R\ \big|\ \kdim\,R/\p\le n\,\big\}$.\\
\item[{\rm (10)}]
$\big\{\,M\in\mod\,R\ \big|\ \ell(M)<\infty\,\big\}\,\leftrightarrow\,\Max\,R$.
\end{enumerate}
\end{ex}

\begin{pf}
In each correspondence, we denote by $\M$ the left-hand subcategory of $\mod\,R$, and by $S$ the right-hand subset of $\Spec\,R$.
Note that it is enough to check either that $\Psi(\M)=S$ or that $\Phi(S)=\M$ since $\Phi$ is an isomorphism with the inverse homomorphism $\Psi$.

(1) The subcategory $\M$ consists of all finitely generated $R$-modules $M$ with $\Ass\,M\subseteq V(I)$, which coincides with $\Phi(V(I))$.
Hence $\Phi(S)=\M$.

(2) Let $M$ be a finitely generated $R$-module.
Note that $\Ass\,M\subseteq\Spec\,R\setminus\Supp\,X$ if and only if $\Ass\,M\cap\Supp\,X=\emptyset$, if and only if $\Ass\,\Hom(X,M)=\emptyset$ (cf. \cite[Exercise 1.2.27]{BH}), if and only if $\Hom(X,M)=0$.
Thus we have $\Phi(S)=\M$.

(3) The equality is well-known.
Putting $X=R/I$ in (2), we obtain the correspondence.

(4) Let $M$ be a finitely generated $R$-module with $\grade(M,X)\ge n$, and let $\p\in\Ass\,M$.
Then there is an $X$-regular sequence $\aa=a_1,\dots,a_n$ in $\Ann\,M$, and $\p\in\Supp\,M$.
Hence $\aa$ is an $X$-regular sequence in $\p$, and we have $\grade(\p,X)\ge n$.
Therefore $\Psi(\M)$ is contained in $S$.
Conversely, if $\p$ is a prime ideal with $\grade(\p,X)\ge n$, then $R/\p\in\M$ and $\p\in\Ass\,R/\p$.
Hence $\p$ is in $\Psi(\M)$.
Therefore $S$ is contained in $\Psi(\M)$, and thus we get $\Psi(\M)=S$.

(5) Let $M$ be an $R$-module.
We have that $\rank\,M=0$ if and only if $M_\p=0$ for every $\p\in\Ass\,R$, if and only if $\Supp\,M\cap\Ass\,R=\emptyset$, if and only if $\Ass\,\Hom(M,R)=\emptyset$, if and only if $\Hom(M,R)=0$, namely $\grade\,M>0$.
Thus the equality holds.
For the correspondence, put $X=R$ and $n=1$ in (4).

(6) Let $\p$ be a prime ideal in $\Psi(\M)$.
Then there exists an $R$-module $M\in\M$ of which $\p$ is an associated prime.
Assume that $\grade(\p,X)>0$.
Then there is an $X$-regular element $a\in\p$, and this is also $M$-regular.
This is a contradiction since $\p\in\Ass\,M$.
Thus $\grade(\p,X)=0$, namely $\p$ belongs to $S$.
On the contrary, let $\p$ be a prime ideal with $\grade(\p,X)=0$.
Then there exists an associated prime $\q$ of $X$ which contains $\p$.
Let $a$ be an $X$-regular element.
Then $a$ is not in $\q$, so is not in $\p$.
Hence $a$ is an $R/\p$-regular element, and $R/\p$ belongs to $\M$.
Since $\p\in\Ass\,R/\p$, the prime ideal $\p$ is in $\Psi(\M)$, and it holds that $\Psi(\M)=S$.

(7) Put $X=R$ in (6), and we get this correspondence.

(8) If $M$ is a finitely generated $R$-module with $\height\,\Ann\,M\ge n$, then $\height\,\p\ge n$ for all $\p\in\Supp\,M$, hence for all $\p\in\Ass\,M$.
Therefore $\Psi(\M)$ is contained in $S$.
If $\p$ is a prime ideal of height at least $n$, then the ideal $\Ann\,R/\p$ of $R$ also has height at least $n$ and $\p\in\Ass\,R/\p$.
Hence $\Psi(\M)$ contains $S$, and thus $\Psi(\M)=S$.

(9) Let $M$ be a finitely generated $R$-module.
Then $\Ass\,M$ is contained in $S$ if and only if $\kdim\,R/\p\le n$ for every $\p\in\Ass\,M$, if and only if $\kdim\,M\le n$.
Therefore $\Phi(S)=\M$.

(10) Putting $n=0$ in (9) yields this correspondence.
\qed
\end{pf}

Note that in the correspondences (1), (4), (5), (8), (9) and (10) in the above example, the left-hand subcategories of $\mod\,R$ are Serre subcategories and the right-hand subsets of $\Spec\,R$ are closed under specializations, hence those correspondences are in fact obtained by the isomorphisms $\phi$ and $\psi$.

\begin{ac}
The author would like to give his deep gratitude to Luchezar L. Avramov and Srikanth Iyengar for valuable discussions and useful suggestions, and thanks Henning Krause and Osamu Iyama very much for helpful comments.
Most of the results stated in this paper were obtained while the author visited the University of Nebraska-Lincoln in the summer of 2006.
The author is also greatly indebted to the University of Nebraska-Lincoln for the hospitality and support.
The author also thanks an anonymous referee for important comments.
\end{ac}


\end{document}